\documentclass[10pt]{article}
\usepackage{amsmath,amssymb,enumerate} %mathtools

\usepackage{amsthm}

\usepackage{tikz}
\usetikzlibrary{arrows}
\usetikzlibrary{positioning}

\newtheorem{theorem}{Theorem}
\newtheorem{lemma}{Lemma}
\newtheorem{proposition}{Proposition}

\newtheorem{remark}{Remark}

\def\ZZ{\mathbb{Z}}
\def\ZZ1{\mathbf{Z}_{\geq 1}}

\def\RR{\mathbb{R}}

\title{Metric and ultrametric inequalities for resistances in directed graphs.
}
\author{
Vladimir Gurvich\thanks{
National Research University Higher School of Economics (HSE) Moscow Russia;
\newline
e-mail: vgurvich@hse.ru and vladimir.gurvich@gmail.com}
}

\begin{document}

\maketitle

\abstract{
Consider an electrical circuit  $G$
each directed edge $e$ of which is a semiconductor
with a monomial conductance function
$y_e^* = f_e(y_e) = y_e^s / \mu_e^r$  if  $y_e \geq 0$  and
$y_e^* =  0$  if  $y_e \leq 0$.
Here  % $e$  is a directed edge,  
$y_e$ is the potential difference (voltage), 
$y_e^*$  is  the current in  $e$, and 
$\mu_e$ is the resistance of $e$; 
furthermore, $r$ and $s$ are 
two strictly positive real parameters common for all edges.
In particular, case $r = s = 1$  corresponds to the Ohm law,
while  $r = \frac{1}{2}, s =1$  may be interpreted as
the square law of resistance typical for hydraulics and gas dynamics.
\newline
We will show that for every ordered pair of nodes $a, b$ of the circuit,
the effective resistance $\mu_{a,b}$  is well-defined.
%and called the resistance distance from  $a$ to $b$.
In other words, any two-pole network with poles  $a$ and $b$
can be effectively replaced by two oppositely directed edges,
from  $a$ to $b$  of resistance $\mu_{a,b}$  and
from  $b$ to $a$  of resistance $\mu_{b,a}$.
\newline
Furthermore, for every three nodes $a, b, c$ the  inequality
$\mu_{a,c}^{s/r} + \mu_{c,b}^{s/r} \geq \mu_{a,b}^{s/r}$  holds,
in which the equality is achieved if and only if
every directed path from  $a$ to $b$  contains  $c$.
\newline
Some limit values of parameters  $s$ and $r$
correspond to classic triangle inequalities. Namely, 
\newline
(i) the length/time  of a shortest directed path,
\newline 
(ii) the inverse width of a bottleneck path, and
\newline
(iii) the inverse capacity (maximum flow per unit time)
\newline
between any ordered pair of terminals  $a$ and  $b$  are assigned to:
\newline
(i) $r = s \rightarrow \infty$,
(ii)  $r = 1, s \rightarrow \infty$,
(iii) $r \rightarrow 0, s = 1$, respectively.
\newline
These results generalize ones obtained in 1987
for  the isotropic monomial circuits,  
% These circuits correspond to symmetric directed graphs,  
% in which pairs of oppositely directed edges have equal resistances
modelled by undirected graphs.  
% (which allows to replace directed graphs by non-directed).
In this special case resistance distances form a metric space,
while in general only a quasi-metric one: symmetry,
$\mu_{a,b} = \mu_{b,a}$  is lost.
\newline
In linear symmetric case
these results are known from 1960-s and  
were generalized to linear the non-symmetric case in 2016-th.   
\newline
MSC classes:
11J83, %Metric theory
90C25, %Convex Programming
94C15, %Applications of graph theory [See also 05Cxx, 68R10]
94C99  %Circuits, Networks: None of the above, but in this section
}

\section{Introduction}
\subsection*{Two-Pole Circuits}
We consider a circuit modeled by a directed graph (digraph) $G = (V,E)$
in which each directed edge $e \in E$
is a semiconductor with the monomial conductivity law
$$y_e^* = f_e(y_e) = y_e^r / \mu_e^s$$
if  $y_e \geq 0$  and
$y_e^* = 0$  if  $y_e \leq 0$.  
Here  $y_e$ is the voltage, or potential difference,
$y_e^* \geq 0$ current, and $\mu_e$  is the resistance of $e$,
while $r$ and $s$ are two strictly positive real
parameters, the same for  all  $e \in E$.

In particular, the case $r = 1$ corresponds to Ohm's law,
while $r = \frac{1}{2}$  is
the so-called  square law of resistance typical for hydraulics and gas dynamics.
In the first  case
$y_e$  is the drop  of potential (voltage) and $y_e^*$  is the current;
in the second case
$y_e$  is the drop of pressure and $y_e^*$  is the flow.
Parameter $s$, in contrast to r, looks redundant, yet,
it plays an important role helping to interpret some limit cases.

Given a circuit  $G = (V,E)$,
let us fix an ordered pair of nodes $a, b \in V$.
We will show that the obtained two-pole circuit $(G, a, b)$
satisfies the same monomial conductivity law.
Let  $y^*_{a,b}$  denote the total current
that comes from  $a$  into  $b$
and $y_{a,b}$ the drop  of potential (voltage) between $a$ and $b$.
It  will be  shown that
$$y^*_{a,b} = f_{a,b}(y_{a,b}) = y^r_{a,b} / \mu^s_{a,b}$$
when  $y_{a,b} \geq 0$,  and
there exists a directed path from  $a$  to  $b$  in  $G$.
If there is no such path  then
$y^*_{a,b} = 0$  for any  $y_{a,b} \geq 0$; in this case we set $\mu_{a,b} = +\infty$.
Also, by convention, we set  $y^*_{a,b} = 0$  when  $y_{a,b} < 0$  or  $a = b$.
In the  latter case $y_{a,b} = 0$ always holds and we set  $\mu_{a,b} = 0$, by convention.

In other words, each two-pole circuit  $(G, a, b)$   can be effectively
replaced by two oppositely directed edges:
from  $a$ to $b$  of resistance $\mu_{a,b}$  and
from  $b$ to $a$  of resistance $\mu_{b,a}$.
Both numbers are  $0$  when  $a = b$.

\subsection*{Main inequality}
For arbitrary three nodes  $a,b,c \in G$,
we will prove
% the following {\em triangle}
inequality

\begin{equation}
\label{main}
\mu^{s/r}_{a, b} \leq \mu^{s/r}_{a, c} + \mu^{s/r}_{c, b}
%%% \;\;\; \forall \;  a,b,c \in V,.
\end{equation}

Furthermore, the inequality in (\ref{main}) is strict
if and only if there exists a  directed path from  $a$  to $b$
that does not contain  $c$.

Obviously, equality holds
if at least one of the three considered resistances
equals  $0$  or  $+ \infty$, that is,
if at least two of the considered three nodes coincide,
or at least one of three directed paths,
from  $a$ to $c$,
from  $c$ to $b$, or
from  $a$ to $b$  fails to exist.
In the latter case, $+ \infty = + \infty$  by convention.

Clearly, if  $s \geq r$  then (\ref{main}) implies the standard metric inequality:

\begin{equation}
\label{triangle}
\mu_{a, b} \leq \mu_{a, c} + \mu_{c, b}
%%% \;\;\; \forall \;  a,b,c \in V,.
%%% \;\; \mbox{whenewver} \;\;  s \geq r
\end{equation}

Thus, a circuit can be viewed as a quasi-metric space in which the distance
from  $a$  to $b$  is the effective resistance  $\mu_{a,b}$.
Note that equality  $\mu_{a,b} = \mu_{b,a}$  holds only
in symmetric case, but may fail in general.

\subsection*{Quasi-metric and quasi-ultrametric spaces corresponding
to asymptotics of parameters $r$ and $s$}
{\bf Playing with parameters  $r$ and $s$}, one can get several interesting
(but well-known)  examples.
Let  $r = r(t)$  and  $s = s(t)$  depend on a real parameter  $t$.
Then, these two functions define a curve in the positive quadrant
$r \geq 0, s\geq 0$.
For the next four limit transitions, as  $t \rightarrow \infty$,
for all pairs of poles  $a, b \in V$, the limits
$\mu_{a,b}  = \lim_{t \rightarrow \infty} \mu_{a, b}(t)$  exist
and can be interpreted as follows:

\begin{itemize}
\item{(i)}
Effective resistance of 
an Ohm semiconductor circuit from  pole  $a$  to pole  $b$; 
$\; s(t) = r(t) \equiv 1,\; $  or more generally, 
$\; s(t) \rightarrow 1$  and  $r(t) \rightarrow 1$.
\item{(ii)}
Standard length (travel time or cost) of a shortest route
from terminal  $a$  to terminal $b$  in a circuit of one-way roads; 
$\; s(t) = r(t) \rightarrow \infty$, or
more generally,  $\; s(t) \rightarrow \infty$  and  $s(t)/r(t) \rightarrow 1$.
\item{(iii)}
The inverse width of a widest bottleneck path 
from terminal  $a$  to terminal  $b$ in a circuit of one way-roads; 
$s(t) \rightarrow \infty$  and  $r(t) \equiv 1$, or
more generally, $r(t) \leq const$, or even more generally
$s(t)/r(t) \rightarrow \infty$.
%We conjecture that the same holds whenever
\item{(iv)}
The inverse capacity (maximum flow per unit time)
from terminal  $a$  to terminal   $b$  in a one-way pipeline;   
$\;\; s(t) \equiv 1$  and  $r(t) \rightarrow 0$;
or more generally, $\; s(t) \rightarrow 1,$ while  $r(t) \rightarrow 0$.
\end{itemize}

\medskip

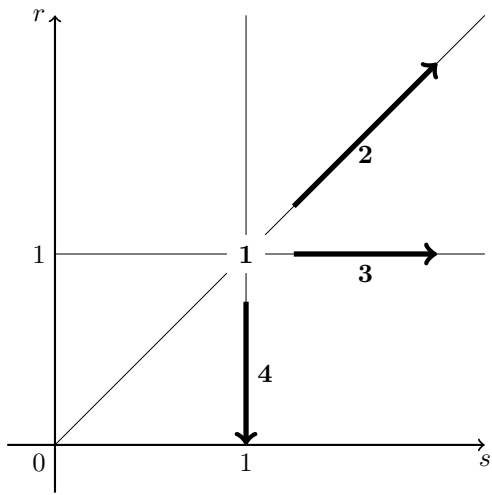
\begin{figure}[h]
\begin{center}
\begin{tikzpicture}
  \draw [->,thick] (-.25in,0in) to (2.25in,0in);
  \draw [->,thick] (0,-0.25in) to (0,2.25in);
  \draw [very thin] (0,1in) to (.9in,1in);
  \draw [very thin] (1.1in,1in) to (2.25in,1in);
  \draw [very thin] (1in,0) to (1in,.9in);
  \draw [very thin] (1in,1.1in) to (1in,2.25in);
  \draw [very thin] (0,0) to (.9in,.9in);
  \draw [very thin] (1.1in,1.1in) to (2.25in,2.25in);
  \draw [->,line width=2pt] (1in,.75in) to node [right] {$\bf{4}$} (1in,0);
  \draw [->,line width=2pt] (1.25in,1in) to node [below] {$\bf{3}$} (2in,1in);
  \draw [->,line width=2pt] (1.25in,1.25in) to node [below] {$\bf{2}$} (2in,2in);
  \node at (1in,1in) {$\bf{1}$};
  \node at (0,1in) [shape=coordinate,label=left:$1$]  {};
  \node at (1in,0) [shape=coordinate,label=below:$1$]  {};
  \node at (0,2.25in) [shape=coordinate,label=left:$r$]  {};
  \node at (2.25in,0) [shape=coordinate,label=below:$s$]  {};
  \node at (0,0) [shape=coordinate,label={below left:$0$}] {};
\end{tikzpicture}

\end{center}
  \caption{Three types of limit transitions for \emph{s} and \emph{r}.}
  \label{pic0}
\end{figure}

All four examples define quasi-metric spaces,
since in all cases  $s(t) \geq r(t)$  for any sufficiently large  $t$
and we assume that  $t  \rightarrow \infty$.
Moreover, for the last two examples the ultrametric inequality

\begin{equation}
\label{ultrametric}
\mu_{a, b} \leq \max(\mu_{a, c}, \mu_{c, b})
%%%  \;\;\; \forall \;  a,b,c \in V,
\end{equation}

\noindent
holds for any three nodes $a,b,c$, because
$s(t)/r(t) \rightarrow \infty$, as
$t \rightarrow \infty$, in cases (iii) and (iv).

\bigskip

These four examples allow us to interpret
$s$  and  $r$  as some important parameters
of transportation problems.

Parameter  $s$  can be viewed as
a measure of divisibility of a transported material;
$s(t) \rightarrow 1$  in examples (i)  and  (iv), because
liquid, gas, or electrical charge are fully divisible;
in contrast,  $s(t) \rightarrow \infty$  for (ii) and (iii), because
a car, a ship, or an individual
travelling from  $a$  to  $b$  is indivisible.

Ratio  $s/r$  can be viewed as
a measure of subadditivity of the transportation cost;
so  $s(t)/r(t) \rightarrow 1$  in examples (i) and (ii), 
because in these cases 
the cost of transportation along a directed path is additive, i.e.,
is the sum of the costs or resistances 
of the directed edges that form this path;
in contrast,  $s(t)/r(t) \rightarrow \infty$  for (iii) and (iv), because
in these cases only edges of the maximum cost
(the width of a bottleneck)  
or capacity of a critical cut matter.

Other values of parameters  $s$  and  $s/r$, between  $1$
and  $\infty$, correspond to an intermediate
divisibility of the transported material and
subadditivity of the transportation cost, respectively.

\bigskip

The metric inequality in case (ii) is obvious.
Let  $\mu_{a,c}$  and  $\mu_{c,b}$  be the lengths
(or the travel times) of the shortest directed paths
$p_{a,c}$  from  $a$  to $c$
and  $p_{c,b}$  from  $c$  to $b$, respectively.
Combining these two paths we obtain a walk  $p_{a,b}$  from  $a$  to $b$.
Thus, the triangle inequality \ref{triangle} follows.

The ultrametric inequality of case (iii) can be proven in a similar way.
Let  $\lambda_{a,c}$  and  $\lambda_{c,b}$  be the largest width
(or wejght) of an object that can be transported
from  $a$  to $c$  and  from  $c$  to $b$, respectively, and
let  $p_{a,c}$) and  $p_{c,b}$) be the corresponding transportation paths.
Combining them we obtain a walk  $p_{a,b}$  from  $a$  to $b$.
Obviously, an object of width   $w$  can be transported along this walk,i.e.,  
$\lambda_{a,b} \geq w$,  whenever 
it  can be  transported from  $a$ to   $c$  and  from $c$ to $b$,
i.e.,  $\lambda_{a,c} \geq w$  and  $\lambda_{c,b} \geq w$.
This implies the  ultrametric inequality (\ref{ultrametric}).

Note that in both above (``indivisible") cases
the inequalities are strict whenever directed paths
$p_{a,c}$  and  $p_{c,b}$  intersect not only in $c$.
Yet, we will show that they may be strict in some other cases too.

``Divisible" case (iv) requires a different approach.
Given a circuit with three fixed nodes  $a,c$, and  $b$,
let  $\lambda_{a,c}$  and  $\lambda_{c,b}$  denote the capacities
(that is, the maximum feasible flows)
from  $a$  to $c$  and from  $c$  to $b$, respectively.
Furthermore, let  $w$  be the  value of a flow
feasible in both cases, in other words,
let inequalities  $\lambda_{a,c} \geq w$  and  $\lambda_{c,b} \geq w$  hold.
Then,  $\lambda_{a,b} \geq w$  holds too, implying (\ref{ultrametric}).
However, to show this, one cannot just combine
(sum up) two flows realizing   $\lambda_{a,c}$  and  $\lambda_{c,b}$,
because  the resulting flow may 
exceed capacities of some edges, thus becoming not feasible.

Instead, inequality
$\lambda_{a,b} \geq w$  can be easily
derived from the classic ``Max Flow - Min Cut Theorem" \cite{FF56}.
According to it, each minimum  $(a-b)$-cut $C: V = V_a \cup V_b$
(such that  $a \in V_a, \; b \in V_b$, and $V_a \cap  V_b = \emptyset$)
is of capacity  $\lambda_{a,b}$.
\begin{itemize}
\item[]
If  $c \in V_a$  then  $C$  is a $(c-b)$-cut too and, hence, $\lambda_{c,b} \leq  \lambda_{a,b}$;
\item[]
if  $c \in V_b$  then  $C$  is a  $(a-c)$-cut too and, hence, $\lambda_{a,c} \leq \lambda_{a,b}$.
\end{itemize}

Thus, $\min(\lambda_{a,c},\lambda_{c,b}) \leq \lambda_{a,b}$,
which is equivalent with  (\ref{ultrametric}).

\subsection*{When (\ref{main}) holds with equality}
In this paper we will prove that (\ref{main})  holds
with equality if and only if
each directed path  from  $a$  to  $b$   contains  $c$.
(For the symmetric case this was shown in \cite{GG92,Gur10}.)
The statement holds for any strictly positive real $r$  and  $s$;
in particular, in case (i), when  $r=s=1$.
Yet, for the asymptotic cases (ii, iii, iv)
only  "if part"  holds, while  "only if" one may fail.
Three examples are as follows:
\begin{itemize}
\item{cases (ii) and (iii).}
Let  $G = (V,E)$  be the directed triangle
in which  $V = \{a,b,c\}$  and  $E = \{(a,b),(a,c),(c,b)\}$.
Obviously,  $G$  contains a directed $(a-b)$-path avoiding  $c$
(just the edge $(a,b)$).
Set  $\mu_{(a,c)} = \mu_{(c,b)} = 1$  and  $\mu_{(a,b)} = 3$.
Then, in case (i) we have
$\mu_{a,c} = \mu_{c,b} = 1$  and  $\mu_{a,b} = 2$.
Thus,  (\ref{triangle}) holds with equality:  $1+1=2$.
In case  $(iii)$  we  have
$\mu_{a,c} = \mu_{c,b} = \mu_{a,b} = 1$.
Indeed, edge $(a,b)$  of the width
$\lambda_{a,b} = 1/3$  is useless.
Thus,  (\ref{ultrametric}) holds with equality: $\max(1,1)=1$.
\item{cases (iv).}
Define digraph $G = (V,E)$  by
$V = \{a,b,c,k,\ell\}$  and
$E  = \{(a,k),(k,c),(c,\ell),(\ell,b),(k,\ell)\}$.
Again  $G$  contains a directed $(a-b$-path avoiding  $c$;
it is  given by vertex-sequence  $a,k,c,\ell,b$.
Set  $\mu_e  = 1$  for all  $e  \in E$.
Then again  $\mu_{a,c} = \mu_{c,b} = \mu_{a,b} = 1$,
since edge $(k,\ell)$  of capacity
$\lambda_{k,l} = 1$  is not needed for transportation,
and  (\ref{ultrametric}) holds with equality: $\max(1,1)=1$.
\end{itemize}

\subsection*{Known special cases of the main inequality}
Our main result (\ref{main}) generalizes some well
(or maybe, not so well) known inequalities
obtained earlier for the following special cases.

\medskip

{\bf Symmetric case.}
A  digraph  $G$  is called {{\em symmetric} if its edges
are split into pairs of oppositely directed edges
$e' = (v', v''), e'' = (v'', v')$.
Respectively, a circuit is called {\em symmetric}
if its graph is  symmetric and
$\mu_{e'} = \mu_{e''}$  for each pair  $e', e''$  introduced above.
In this case one can replace each such pair  $e', e''$
by a non-directed edge  $e$,
thus, replacing the digraph of the circuit by a non-directed graph.
For this case, the main inequality (\ref{main}) was shown in \cite{GG87};
see also  \cite{Gur10,Gur12,GV12} for more details.

The equality holds in (\ref{main}) if and only if
every path from  $a$  to  $b$  contains  $c$.
First it was shown in Section 16.9 of \cite{GG92}; see also  \cite{Gur10}).
[It was also shown in Section 16.9  that
the  monomial conductance law is the only only when
the  effective resistance  $\mu_{a,b}$  of the  two-pole
circuit $(G,a,b)$ is a real number.
In general, it is a monotone non-decreasing function
for an arbitrary monotone circuit \cite{Min60}.]

Clearly, equality $\mu_{a,b} = \mu_{b,a}$  holds in the isotropic (symmetric) case.
Thus, resistance distances of symmetric circuits form a metric spaces.
Yet, in anisotropic (non-symmetric) case the above equality  may fail and
we obtain only quasi-metric spaces, in general.

\bigskip

{\bf Linear case}, $r=s=1$.
In the symmetric linear case, the metric resistance inequality
was discovered by Gerald E. Subak-Sharpe  \cite{Sha67,Sha67a};
see also \cite{MS68, Shap87,SS97,Sub89,Sub90,Sub91,Sub92,CS97,CS00,Che11} and
preceding works \cite{You59,SS60,SR61,SS65,SS67}.
This result was rediscovered several times later.

\medskip

Let us notice that the proof of (\ref{main})
given in \cite{GG87,Gur10} for
the isotropic monomial conductance differs a lot from
the proof of \cite{Sha67,Sha67a} for the symmetric linear case.
In this paper we extend the first proof to the anisotropic case
or, in other words, to digraphs.
To  make the presentation self-contained we copy here
some parts of \cite{Gur10}.

\medskip

For the linear non-symmetric case
the quasi-metric inequality, along with many related results,
was recently obtained in \cite{YSL16}.

\subsection*{Continuum.}  It would be natural to conjecture that
the above approach can be developed
not only for the discrete  circuits but for continuum as well:
inequality (\ref{main}) and its corollaries
should hold in this case  too.
Sooner or later, this will become the subject of a separate research.
The same four subcases appear:
linear and monomial, isotropic and anisotropic.

\section{Resistances of two-pole circuits}
\label{s1}
% ; proof of inequality
% $\mu^{s/r}_{a, b} \leq \mu^{s/r}_{a, c} + \mu^{s/r}_{c, b}$}

\subsection*{Conductance law}
Let  $e$  be a semi-conductor with the monomial conductivity law

\begin{equation}
\label{conduct}
y_e^* = f_e(y_e) = \lambda_e^s y_e^r = \frac{y_e^r} {\mu_e^s}
\mbox{\;\;if\;\;} y_e \geq 0 \mbox{\;\; and \;\;} 0 \mbox{\;\; if  \;\;} y_e \leq 0.
\end{equation}

\noindent
Here $y_e$  is the {\em voltage or potential difference},
$y_e^*$  {\em current},
$\lambda_e$  {\em conductance}, and
$\mu_e = \lambda_e^{-1}$  {\em resistance} of  $e$;
furthermore,  $r$  and  $s$  are two
strictly positive real parameters independent of  $e$.
Obviously, the monomial function  $f_e$  is continuous,
strictly monotone increasing when $y_e \geq 0$,
and taking all non-negative real values.

%\begin{itemize}
%\item{} continuous, strictly monotone increasing, and taking all real values;
%\item{} symmetric (odd or isotropic), that is,  $f_e(-y_e) = -f_e(y_e)$;
%\item{} the inverse function  $f_e^{-1}$  is also monomial
%with parameters  $r' = r^{-1}$  and  $s' = s^{-1}$.
%\end{itemize}

\subsection*{Main variables and related equations}
\label{s2}
A semi-conductor circuit is modeled by
a weighted digraph  $G = (V, E, \mu)$
in which weights of the edges are
their {\em positive} resistances  $\mu_e, e \in E$.

%%  and by two positive real parameters
%$r = r(G)$ and  $s = s(G)$  independent of $e$.

Let us introduce the following four groups of real variables;
two for each $v \in V$  and $e \in E$:
{\em potential}  $x_v; \,$  {\em difference of potentials, or voltage} $y_e; \,$
{\em current}  $y_e^*; \,$  {\em sum of currents, or flux}  $x_v^*$.

\medskip

The above variables are not independent.
By (\ref{conduct}), current  $y_e^*$  depends on voltage  $y_e$.
Furthermore, the voltage (respectively, flux) is
a liner function of the potentials (respectively, of the currents).
These functions are defined by the node-edge incidence function of the digraph  $G$:

\begin{equation}
  \mbox{inc}(v,e)=\left\{
  \begin{aligned}
    +1&, &&\text{if node} \;\; v\ \;\; \text{is the beginning of} \;\; e;\\
    -1&, &&\text{if node} \;\; v\ \;\; \text{is the end of} \;\; e;\\
     0&, &&\text{if} \;\; v \;; \text{and} \;\;  e \text{are not incident}.
  \end{aligned}\right.
\end{equation}

We will assume that the next two systems of linear equations always hold:

\begin{gather}
  y_{e}=\sum_{v\in V} \mbox{inc}(v,e)x_v;
  \label{eq1}\\
  x^*_v=\sum_{e\in E} \mbox{inc}(v,e)y_e^*.
  \label{eq2}
\end{gather}

Let us notice that equation (\ref{eq1})
for a directed edge  $e = (v', v'')$  can be reduced to
% $y_e= \mbox{inc}(e, v') x_{v'}+ \mbox{inc}(e, v'') x_{v''}$
% and even further to
$y_e = x_{v'} - x_{v''}$.
% yet, for the latter it should be assumed that
% $e$  is directed from  $v'$  to  $v''$.

We say that the first Kirchhoff law holds
for a node  $v$  if  $x_v^* = 0$.

\medskip

Let us introduce four vectors, one for each group of variables:

$$x = (x_v  \mid  v \in V), \; x^* = (x_v^*  \mid  v \in V), \;\;
y = (y_e \mid e \in E), \; y^* = (y_e^* \mid e \in E),$$
$$x, x^* \in \RR^n; y, y^* \in \RR^m,$$

\noindent
where $n = |V|$  and   $m = |E|$  are the numbers of nodes and edges
of the digraph $G = (V, E)$.
Let   $A = A_G$  be the edge-node  $m \times n$
incidence matrix of graph  $G$, that is,
$A(v, e) = inc(v, e)$  for all $v \in V$  and  $e \in E$.
Equations (\ref{eq1}) and  (\ref{eq2}) can be rewritten
in this matrix notation as
$y = Ax$  and  $x^* = A^T y^*$, respectively.

It is both obvious and well known that
these two equations imply the following chain of identities:

\medskip
% \begin{equation}
$$(x, x^*) = \sum_{v\in V}x_v x_v^*=\sum_{e\in E}y_e y_e^* = (y, y^*).$$
% \end{equation}

Recall that  $y*$  is uniquely defined by   $y$
according to the conductance law (\ref{conduct}).
Thus, given vector $x$, the remaining three vectors  $y$, $y^*$, and  $x^*$
are uniquely defined by $x$
(\ref{eq1},\ref{conduct},\ref{eq2}).
This triple will not change
if we add an arbitrary real constant $c$  to all coordinates  of   $x$, while
multiplying  $x$  by  $c$  will result in
multiplying  $y$  by  $c$  and  $y^*, x^*$  by  $c^r$.
More precisely, the following scaling property cllearly holds.

\begin{lemma}
\label{l0}
For any positive constant  $c$, two quadruples
$(x, y, y^*, x^*)$   and
$(cx, cy, c^r y^*, c^r x^*)$  can satisfy
all equations of (\ref{eq1},\ref{eq2},\ref{conduct})
only simultaneously.       \qed
\end{lemma}

\subsection*{Two-pole boundary  conditions}
In general theory of monotone circuits,
one can consider arbitrary monotone functions:
a non-decreasing one  $y^*_e =  f_e(y_e)$
for each $e \in E$  and
a non-increasing one  $x^*_v = g_v(x_v)$
for each $v \in V$; see
\cite{Duf47,Min60,Roc67,Roc70} and also
\cite{GG90,GG92,OGG86a,OGG86b,OGG90}.

In case of the two-pole circuits
we restrict ourselves by the monomial conductance  law  (\ref{conduct}).
Fix an {\em ordered} pair of poles  $a,b$,  potentials
\begin{equation}
\label{boundary}
x_a = x_{a}^0, x_b = x_b^0,
\end{equation}
in them, and require the first Kirchhoff law for any other node:

\begin{equation}
\label{Kirchhoff}
x^*_v = 0, \;\; v \in V \setminus \{a,b\}.
\end{equation}

By convention,  $y^*_{a,b} = 0$  if  $a = b$   or   $y_{a,b} = x_a - x_b  \leq  0$.
So, w.l.o.g. we can assume that  $a \neq  b$  and  $x_{a}^0 \geq  x_{b}^0$.

\begin{remark}
\label{r0}
By Lemma  \ref{l0}, it would be sufficient
to replace (\ref{boundary}) by
$x_{a}^0 = 1$  and  $x_{b}^0 = 0$.
It would be also possible to replace it by
$x^*_a = x^{*0}_a$  (or  $x^*_a = 1$).
Then,  $x^*_b = - x^{*0}_a$
(resp., $x^*_b = -1$) will automatically hold, by (\ref{Kirchhoff}).
\end{remark}

We call  a vector $x = x(a,b)$  a solution
of the two-pole circuit  $(G, a, b)$  if
the corresponding quadruple $(x, y, y^*, x^*)$
satisfies all equations (\ref{conduct} - \ref{Kirchhoff}).

\bigskip

In \cite{Gur10} the monomial symmetric case was considered and
it was shown that
there exists a unique solution  $x = x(a,b)$
whenever  $a$  and  $b$  belong to the same connected
component of  $G$.

Yet, this claim cannot be extended directly to digraphs.
For example, let  $(G, a, b)$  be a  directed $(a,b)$-$k$-path
$a = v_0, v_1, \ldots, v_k = b$   from  $a$  to $b$.
Then potential $x$  is unique if  $x^0_a \geq x^0_b$, but otherwise,
when  $x^0_a < x^0_b$, any non-decreasing
$x^0_a = x_0 \leq x_{v_1} \leq  \ldots \leq x_{v_k} = x^0_b$
will be a solution, with no current, that is,
$y_e^* = 0$  for all  $e = (v_{j-1}, v_j)$   for  $j = 1, \ldots, k$.

In the directed case we will prove that $y^*$
(rather than  $x$)  is the same in all solutions.
Furthermore, let $G^+$  be the subgraph of $G$  defined by
all directed edges  $e \in E$  such that  $y^*_e > 0$.
Then in all solutions potentials  $x$
are uniquely defined on vertices of $G'$.

\subsection*{Existence of a solution}
We will apply Method of Successive Approximation (MSA)
increasing potentials of some nodes, one by one in a certain order.

Obviously, when we increase  $x_v$
(keeping all remaining potentials  $x_u$  unchanged)
the corresponding flux  $x^*_v$  is non-decreasing;
furthermore, it is strictly increasing if and only if
$G$  contains an edge $(v, u)$  with $x_v \geq x_u$  or
an edge $(u, v)$  with $x_v \leq x_u$.
Respectively, $x^*_u$  in any other node
$u \in V \setminus \{v\}$   is non-increasing;
furthermore, it is strictly decreasing if and only if
$(v, u)$  is an edge and  $x_v \geq x_u$  or
$(u, v)$  is an edge and  $x_v \leq x_u$.

\medskip

Let us set  $x_a = x_a^0$ and $x_v = x_b^0$
for all nodes  $v \in V \setminus \{a\}$, including  $b$.
In the course of iterations, potentials
$x_a = x_a^0$ and $x_b = x_b^0$  will remain unchanged,
while  the all other potentials  $x_v$, on the nodes from
$W = V \setminus \{a,b\}$,  will be recomputed by MSA as follows.

Order arbitrarily the nodes of  $W$ and
consider them one by one in this order repeating cyclically.
If  $x^*_v = 0$,  skip this node and go to the next one.
If  $x^*_v < 0$, increase  $x_v$  until $x^*_v$ becomes  $0$.
The latter is possible, since, by (\ref{conduct}),
$f_e$  is continuous and
$y^*_e \rightarrow + \infty$  as  $y_e \rightarrow + \infty$.
Let us notice that $x^*_v$ may remain  $0$   for some time, but
we stop increasing  $x_v$  the first moment when  $x^*_v$  becomes  $0$,
and proceed to the next node.

The following claims can be easily
proven together by induction on the number of iterations.

\begin{itemize}
\item{(i)}
For any node $v \in W$  its potential $x_v$
is monotone non-decreasing and it remains bounded by  $x_a^0$  from above.
Hence, it tends to a limit  $x_v^0$  between  $x_a^0$  and  $x_b^0$.
\item{(ii)}
These limit potentials solve the two-pole circuit $(G,a,b)$.
\item{(iii)}
Fluxes  $x^*_v$  remain non-positive for all $v \in V \setminus \{a\}$.
In contrast,  $x^*_a$  remains non-negative.
Furthermore, $x^*_a$  and  $x^*_b$   are monotone non-increasing.
\item{(iv)} For the limit values of potentials and fluxes we have:
\newline
If  $G$ contains no directed  path from  $a$  to  $b$
then $x_v^{*0} = 0$  for all $v \in V$, including  $a$ and $b$.
In this case  $x^0_v = x_a^0$  if
$G$  contains a directed path from  $a$  to  $v$,
otherwise  $x^0_v = 0$.
% \newline
If  $G$  contains a directed  path from  $a$  to  $b$
then   $x_a^{*0} > 0$.
(Respectively, $x_b^{*0} = - x_a^{*0} < 0$ and
$x_v^{*0} =  0$  for all  $v \in W$,
in accordance with (\ref{Kirchhoff}.)
\item{(v)}
The limit potentials  $x_v^0$  take
only values $x_a^0$ and $x_b^0$  if and only if
$G$ contains no directed  path from  $a$  to  $b$  or
every such path consists of only one edge.
\end{itemize}

Existence of a solution is implied by  (ii).

\begin{remark}
\label{pumping}
A very similar monotone potential reduction
(pumping) algorithm
for stochastic games with perfect information
was suggested in \cite{BEGM10,BEGM15}.
\end{remark}

\subsection*{Uniqueness of the solution}
A solution  $x$  of a two-pole  circuit $(G,a,b)$  may be not unique.
Suppose that $G$ contains an induced directed path
$P$  from  $u$ to $v$  of length greater than  $1$  and that
$(G,a,b)$   has a solution  $x$  with  $x_u < x_v$.
Then every monotone non-decreasing sequence of potentials on  $P$  is feasible.
Notice, however, that the current along  $P$  will be zero for any such sequence.

We will demonstrate that, in general,
vector of currents $y^*$
(and, hence, $x^*$  too) is unique for all solutions of  $(G,a,b)$.
It follows directly from an old classical result
relating solutions of an arbitrary monotone circuit
with a pair of dual problems of convex programming
\cite{Duf47,Min60,Roc67,Roc70}; see also
\cite{OGG86a,OGG86b,OGG89,OGG90,GG87,GG90,GG92}.

First, note that, by Lemma \ref{l0}},
we can replace the  boundary conditions (\ref{boundary}) by
\begin{equation}
\label{boundary*}
x^*_a = x^{*0}_a, \;\; x^*_b = - x^{*0}_a,
\end{equation}
and recall that   $x_v^* = 0$  for all  $v \in W = V \setminus \{a,b\}$  by (\ref{Kirchhoff}).

The Joule-Lenz heat on  $e$  is defined by the current  $y_e^* \geq 0$
as the integral

\begin{equation}
\label{dissipation-e}
F^*_e(y^*_e) = \int f^{-1}_e (y^*_e) \, d y^*_e  =
\frac {\mu_e^{s/r}} {1 + 1/r} y_e^{* \; {1 + 1/r}},
\end{equation}

\noindent
which is a strictly convex function of  $y_e^* \geq 0$.
Furthermore, the total heat dissipated in the circuit is additive:
\begin{equation}
\label{dissipation}
F^*(y^*) = \sum_{e \in E} F^*_e (y^*_e).
\end{equation}

\noindent
It is a strictly convex
function of  $y^*$  defined on the positive ortant, $y^* \geq 0$.

By the classical results
(see, for exampple, Rockafellar  \cite{Roc67})
solving  $(G,a,b)$  is equivalent with minimizing
dissipation  $F^*(y^*)$   subject to the following constraints

\begin{equation}
\label{boundary*}
x^*_a = x^{*0}_a, \; (x^*_b = - x^{*0}_a),  \; x^*_v = 0
\; \forall \; v \in V \setminus \{a,b\}, \;\; and \;\;  y^*  \geq 0.
\end{equation}

Since  $x^* = A y^*$, we obtain the minimization problem
for a strictly convex function of  $y^*$
subject to linear constraints  on  $y^*$.
It is known from calculus that solution $y^{*0}$  is unique in this case.

\begin{remark}
In the non-directed (isotropic) case
the above equivalence is well-known in physics as
the minimum dissipation principle.
It is applicable for arbitrary ``bounddary conditions"
not only to the two-pole circuits.
\end{remark}

Thus, all solutions of  $(G,a,b)$
have the same current vector  $y^{*0}$.
This implies the uniqueness of the flux vector
$x^{*0} = A^* y^{*0}$  as well.

\medskip

Let us denote by  $G^+ = (V^+, E^+)$  the subgraph of  $G$
formed by the edges  $e \in E$  with  positive currents, $y^{*0} > 0$.
The following properties of   $G^+$  are  obvious:

\begin{itemize}
\item{(j)} Digraph $G$  contains a directed path from  $a$  to  $b$
if and only if  $G^+$  is not empty.
\item{(jj)} In the latter case it contains the poles, $a,b \in V^+$,
and at least one directed path from  $a$ to $b$,
but not necessarily all such pathes.
Yet, any vertex or edge of  $G^+$  belongs to such a path.
% \item{(jj)} Potentials on  $V^+$ are uniquely defined,
% up adding a constant to all of them .
% Kakaja konstanta?!  Est' zhe granichnyye uslovija
\item{(jjj)}  Potentials are strictly decreasing on each edge of  $G^+$
and, hence, it has no directed cycles.
\end{itemize}

\subsection*{Conductance functions of two-pole networks}
Given a two-pole network  $(G,a,b)$, let us define
the potential drop and current from  $a$  to  $b$  as
$$y_{a,b} = x_a - x_b, \;\;  y^*_{a,b} = x^*_a = - x^*_b.$$

If there is no directed path from  $a$  to  $b$  in  $G$,
let us set  $\mu_{a,b} = + \infty$,
since in this case  $y^*_{a,b} = 0$  for any  $y_{a,b}$.
Otherwise, by Lemma \ref{l0}, $y^*_{a,b}$ depends on $y_{a,b}$
as in  (\ref{conduct}):

\begin{equation}
\label{conduct-ab}
y_{a,b}^* = f_{a,b}(y_{a,b}) = \lambda_{a,b}^s y_{a,b}^r = \frac{y_{a,b}^r} {\mu_{a,b}^s}
\mbox{\;if\;} y_{a,b} \geq 0 \mbox{\; and \;} 0 \mbox{\; if  \;} y_{a,b} \leq 0.
\end{equation}

Two strictly positive real values
$\lambda_{a,b}$  and  $\mu_{a,b} = \lambda_{a,b}^{-1}$ are called
{\em conductance} and, respectively, {\em resistance} of  $(G,a,b)$.

\begin{remark} It is shown in Section 6.9 of \cite{GG92} that
among all monotone conductance laws
the monomial one  is the only case when resistance of
a two-pole network is a real number; in other words,
up to a real factor, 
the same function  $f$  describes the conductances
$f_e(y_e)$  and  $f_{a,b}(y_{a,b})$.
\end{remark}

\subsection*{Monotonicity of effective resistances and Braess' Paradox}
\label{mon}
Given a two-pole circuit  $(G, a, b)$, where  $G = (V, E, \mu)$,
let us fix an edge  $e_0 \in E$, replace the resistance  $\mu_{e_0}$
by a smaller one,  $\mu'_{e_0} \leq \mu_{e_0}$, and denote by  
$G' = (V, E, \mu')$  the obtained circuit.

Of course, the total resistance will not increase either,
that is, $\mu'_{a,b} \leq \mu_{a,b} $  will hold.
Yet, how to prove this  "intuitively obvious" statement?
Somewhat surprisingly, the simplest way is 
to apply the the minimum dissipation principle again; 
see, for example, \cite{Lya85},

Let  $y^*$  and  $y^{*'}$  be the (unique) current vectors that
solve  $(G, a, b)$  and  $(G', a, b)$, respectively.
Since $\mu'_{e_0} \leq \mu_{e_0}$, inequality
$F^{*'}_{e_0}(y^*_{e_0}) \leq  F^*_{e_0}(y^*_{e_0})$  is implied by
(\ref{dissipation-e}).
Furthermore, $F^{*'}_e(y_e) = F^*_e(y_e)$  for all other  $e \in E$,
distinct from  $e_0$.
Hence,  $F^{*'}(y^*) \leq  F^*(y^*)$ holds by (\ref{dissipation}).
As we know, all solutions of   $(G', a, b)$
have the same vector of currents  $y^{*'}$,
which may differ from  $(y^*)$ and,
by the minimum dissipation principle,
% the dissipation in the first case
% is at  most the dissipation in the second one:
we have  $F^{*'}(y^{*'}) \leq F^{*'}(y^*)$.
From this, by transitivity, we conclude that
$F^{*'}(y^{*'}) \leq F^*(y^*)$  and,
by  (\ref{dissipation-e},\ref{dissipation}),
conclude that  $\mu'_{a,b} \leq \mu_{a,b}$  holds.

In particular, when an edge  $e$  is eliminated from  $G$,
its finite resistance  $\mu_e$ is replaced by $\mu'_e = +\infty$.
It was just shown that, by this operation,
the effective resistance is not reduced, that is,
$\mu_{a,b} \leq \mu'_{a,b}$  holds.

\medskip

In general, for monotone circuits \cite{Duf47,Min60}
the conductance function of its edge
may be an arbitrary, not necessarily  monomial, monotone non-decreasing function:

\smallskip

$y_e^* = f_e(y_e)  % \lambda_e^s y_e^r = \frac{y_e^r} {\mu_e^s}
\mbox{\;\;if\;\;} y_e \geq 0 \mbox{\;\; and \;\;} 0 \mbox{\;\; if  \;\;} y_e \leq 0,$

\smallskip

Then, by results of  \cite{Min60}, the conductance law of
a two-pole network $(G,a,b)$  is represented by a similar formula:

\smallskip

$x^*_a = - x^*_b = y_{a,b}^* = f_{a,b}(y_{a,b})  % \lambda_e^s y_e^r = \frac{y_e^r} {\mu_e^s}
\mbox{\;\;if\;\;} y_{a,b} \geq 0 \mbox{\;\; and \;\;} 0 \mbox{\;\; if  \;\;} y_{a,b} \leq 0$,

\smallskip
\noindent
where  $f_{a,b}$  is a monotone non-decreasing function too.

When we reduce conductance function  $f_e(y_e)$  of an edge  $e \in E$,
the effective conductance function  $f_{a,b}$  may increase for {\em some}
(certainly, not for all) values of its  argument  $y_{a,b} = x_a - x_b$.
This phenomenon is known as Braess paradox \cite{Bra69}.

The above monotonicity principle implies that this paradox is not possible
for circuits with the monomial conductance law
provided parameters  $r$  and  $s$  are the same for all edges $e \in E$.
Indeed, in this case resistance $\mu_{a,b}$  between the poles
is a real number; moreover, 
it is a monotone function of resistances $\mu_e$  of edges $e \in E$, 
as it was shown above.
Yet, the  paradox can appear for monomial circuits
in which parameter $r = r_e$  depends on $e$, or when
some edges have non-monomial monotone conductance functions.

\section{Proof of the main inequality and some related claims}
\label{s3}

Here we prove our main result
generalizing the triangle inequality of \cite{GG87,Gur10}
from graphs to digraphs as follows.

\begin{theorem}
\label{t0}
Given a weighted digraph  $G = (V, E, \mu)$
with strictly positive weights-resistances $(\mu_e \, | \, e \in E)$,
three arbitrary nodes  $a,b,c \in V$, and
strictly positive real parameters  $r$  and  $s$,  
inequality (\ref{main}) holds:
$\;\; \mu^{s/r}_{a, b} \leq \mu^{s/r}_{a, c} + \mu^{s/r}_{c, b}$.
% \newline
Moreover, it holds with equality if and only if
node $c$  belongs to
{\em every} directed path from  $a$  to  $b$  in  $G$.
\end{theorem}

%\begin{remark}
%The proof of the first statement was sketched in \cite{GG87};
%see also \cite{Gur10}. Both claims were proven in \cite{GG92}.
%Here we shall follow the plan suggested in \cite{GG87} but give
%more details.
%\end{remark}

\proof
W.l.o.g. we can assume that
$G$  contains directed paths
from   $a$  to  $c$   and  from  $c$ to  $b$.
Indeed, otherwise
$\mu_{a,c}$  or  $\mu_{c,b}$  is  +$\infty$
and there is nothing to prove.
By this assumption,
$G$  contains a directed walk from  $a$  to $b$
passing through  $c$.
Hence, $G$  also contains a directed path  from
$a$ to $b$, but the latter may avoid  $c$.

Anyway, there exists a
(not necessarily unique) solution of $(G,a,b)$
for any fixed potentials $x_a^0, x_b^0$ in the poles  $a,b$.
If  $x_a^0 \leq  x_b^0$  then  $y^*_{a.b} = 0$, by definition, and
again there is nothing to prove.
Thus, w.l.o.g. we assume that  $x_a^0 > x_b^0$
(Moreover, we could assume w.l.o.g. that  $x_a^0 = 1$ and  $x_b^0 =0$,
but will not do this.)

We make use of the same arguments as in subsection
"Existence of a solution".
Consider a solution $x^0 = x^0(a,b)$
constructed there and denote by  $x_c^0$  the obtained potential in  $c$.
By construction, $x_a^0 \geq x_c^0 \geq x_b^0$
and at least one of these two inequalities is strict.
(Actually, such inequalities hold for any solution  $x = x(a,b)$,
but one chosen  $x_c^0$  will be enough for our purposes.)

Now let us consider the two-pole circuit  $(G, a, c)$
and fix in it  $x_a = x_a^0$   and  $x_c = x_c^0$,
standardly requiring the first Kirchhoff law,
$x_v^* = 0$  for all other vertices
$v \in W = V \setminus {a,c}$, including  $v = b$.

\begin{lemma}
\label{l3(ac)}
The obtained currents in the circuits $(G, a, b)$ and $(G, a, c)$
satisfy inequality
$\;\;\; y^*_{a,b} \geq  y^*_{a,c}.$
Moreover, the equality holds if and only if
$c$  belongs to every directed path from  $a$  to  $b$.
\end{lemma}

\proof
First, recall that the current vectors
(and, hence, the values of
$x^*(a,b) = y^*_{a,b}$  and  $x^*(a,c) = y^*_{a,c}$ too)
are well-defined, that is, remain the same for any solutions
$x(a,b)$  and  $x(a,c)$  of  $(G, a, b)$ and $(G, a, c)$
with boundary conditions  $x_a^0, x_b^0$  and
$x_a^0,  x_c^0$, respectively.

Again we apply  MSA  to compute  $x(a,c)$, yet,
this time we take  $x(a,b)$  as the initialization.
Thus, in the beginning we have

\smallskip

$x^*_a(a,b) = -x^*_b(a,b), \; x^*_c(a,b)=0$

\smallskip
\noindent
and at the end we will have

\smallskip

$x^*_a(a,c) = -x^*_c(a,c), \; x^*_b(a,c)=0$.

\smallskip

Potentials  $x_a = x^0_a$  and  $x_c = x^0_c$ satisfy the boundary conditions
and will stay unchanged in the course of iterations, while
the remaining potentials  $x_v$  on the nodes from
$W = V \setminus \{a,c\}$
will be determined by MSA as follows.
Order arbitrarily the nodes of  $W$ and
consider them one by one in this order repeating cyclically.
If  $x^*_v = 0$,  skip this node  $v$  and go to the next one.
If  $x^*_v < 0$, increase  $x_v$  until
(the very first moment when)  $x^*_v$ becomes  $0$.
Then proceed with the next node.
The following claims can be easily
proven together, by induction on the number of iterations.

\begin{itemize}
\item{(i)} Fluxes  $x_a^*$  and  $x^*_c$ are both monotone non-decreasing
and remain non-negative and non-positive, respectively.
Moreover, $x_v^*$  remain non-positive
for all $v \in V \setminus \{a\}$.
\item{(ii)} All potentials  $x_v$  are monotone non-decreasing
and remain bounded  by  $x_a^0$  from above.
Hence, $x_v$  tends to a limit  $x_v^0(a,c)$
between  $x_a^0$  and  $x_v^0(a,b)$.
\item{(iii)}  These limit potentials solve the two-pole circuit $(G,a,c)$.
\item{(iv)} Since  $G$  contains directed  pathes
from  $a$  to  $c$  and  from   $c$  to  $b$,
for the limit values of the fluxes we have:
$x_a^{*0}(a,c) > 0, \;  x_c^{*0}(a,c) = - x_a^{*0}(a,c)$, and
$x_v^{*0}(a,c) =  0$  for all  $v \in V \setminus \{a,c\}$,
in accordance with (\ref{Kirchhoff}).
\item{(v)} $x_a^{*0}(a,b) = y^*_{a,b} \geq  y^*_{a,c} = x_a^{*0}(a,c)$
\end{itemize}

Existence of a solution is implied by  (iii).
The last inequality holds, because
potential  $x_a$  is constant, while all other
potentials are not decreasing.
Hence, the  flux from  $a$  cannot increase.

\bigskip

If all directed pathes from  $a$  to  $b$  contain  $c$
then equality holds in  (v).
Indeed, in this case $x_v$  will not be changed by MSA for
any vertex  $v$  that belongs to a path from  $a$  to  $c$.
Hence, the flux  $x^*_a$  remains constant, resulting in
$y^*_{a,b} =  y^*_{a,c}$.

Suppose conversely that
$G$  contains a directed path $P$  from  $a$  to  $b$  avoiding  $c$.
Then, order the nodes of  $W = V \setminus \{a,c\}$  so that
the nodes of  $P$  go first ordered from  $b$  to  $a$.
Obviously, in  $|P| - 1$ steps the flux  $x^*_a$
will be strictly reduced.
Furthermore, $x^*_a$  is non-increasing,
and its initial and limit values are
$y^*_{a,b}$  and   $y^*_{a,c}$, respectively.
Thus, $y^*_{a,b} >  y^*_{a,c}$.
Recall that these two numbers are well-defined,
although solutions of  $(G,a,b)$   and   $(G,a,c)$ 
are not necessarily unique.  This  proves the lemma.
\qed

\begin{remark}
The same arguments prove  that
inequality  

\smallskip

$\; f_{a,b}(x_a -  x_b) = y_{a,b}^* \geq  y_{a,c}^* = f_{a,c}(x_a - x_c)$ 

\smallskip
\noindent
holds not only for monomial but
for arbitrary monotone non-decreasing conductivity functions.
\end{remark}

In the exactly same way we can apply MSA to $(G,c,b)$
again taking a solution of  $(G,a,b)$
as an initial approximation.
Clearly this will result  in inequality
$\;\;\; y^*_{a,b} \geq  y^*_{a,c}$, in which
the equality holds if and only if
$c$  belongs to every path between  $a$   and  $b$.
Summarizing we obtain the following statement:

\begin{proposition}
\label{p-abc}
For an arbitrary weighted digraph  $G$
and nodes $a,b,c$  in it, the inequality
$y^*_{a,b} \geq   max(y^*_{a,c}, y^*_{c,b})$
holds  and the following five statements are equivalent:
\begin{itemize}
\item{(ac)} $y^*_{a,b} = y^*_{a,c}$;
\item{(bc)} $y^*_{a,b} = y^*_{c,b}$;
\item{(ab)} $y^*_{a,c} = y^*_{c,b}$;
\item{(acb)} every  directed path from  $a$  to  $b$  contains  $c$;
\item{(=)} $\mu^{s/r}_{a,b} = \mu^{s/r}_{a,c} + \mu^{s/r}_{c,b}$.
\end{itemize}
\end{proposition}

For the rest of the proof of Theorem \ref{t0}
we will need only elementary "high-school" transformations:

\begin{equation}
\label{currents}
y^*_{a,b} = \frac {(x_a^0 - x_b^0)^r} {\mu_{a, b}^s}
\geq   \frac {(x_a^0 - x_c^0)^r} {\mu_{a, c}^s} = y^*_{a,c};  \;\;\;
\newline
y^*_{a,b} = \frac {(x_a^0 - x_b^0)^r} {\mu_{a, b}^s}
\geq   \frac {(x_c^0 - x_b^0)^r} {\mu_{c, b}^s} = y^*_{c,b},
\end{equation}

\noindent
which can be obviously rewritten as follows

\begin{equation}
\label{currensA}
\left(\frac {\mu_{a, c}} {\mu_{a, b}}\right)^{s/r}  \geq
  \frac {x_a^0 - x_c^0} {x_a^0 - x_b^0}; \;\;\;
\left(\frac {\mu_{c, b}} {\mu_{a, b}}\right)^{s/r}  \geq
  \frac {x_c^0 - x_b^0} {x_a^0 - x_b^0}
\end{equation}

Summing up these two inequalities we obtain (\ref{main}).

\medskip

The  above computations show that
(\ref{main}) holds with equality if and only if
$y_{a,b}^* = y_{a,c}^*$  and  $y_{a,b}^* = y_{c,b}^*$.
By Propositions \ref{p-abc} these two equations
are equicallent and hold
if and only if  $c$  belongs to each directed path from  $a$  to  $b$.
\qed

\section{Three limit cases}
\subsection*{Parallel and series connection of edges}
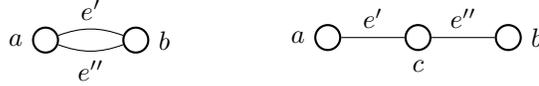
\begin{figure}[h]
\begin{center}
\parbox[c]{1.2in}{
\begin{tikzpicture}[scale=.6]
  \node at (-1,0) [shape=circle,minimum size=6pt,draw=black,thick,label=left:$a$] (a) {};
  \node at (1,0) [shape=circle,minimum size=6pt,draw=black,thick,label=right:$b$] (b) {};
  \draw (a) to [out=20,in=160] node [above] {$e'$} (b);
  \draw (a) to [out=-20,in=-160] node [below] {$e''$} (b);
\end{tikzpicture}
}\qquad
\parbox[c]{2in}{
\begin{tikzpicture}[scale=.6]
  \node at (-2,0) [shape=circle,minimum size=6pt,draw=black,thick,label=left:$a$] (a) {};
  \node at (0,0) [shape=circle,minimum size=6pt,draw=black,thick,label=below:$c$] (c) {};
  \node at (2,0) [shape=circle,minimum size=6pt,draw=black,thick,label=right:$b$] (b) {};
%  \node at (2,0) [label=right:$b$] (b) {};
  \draw (a) to node [above] {$e'$} (c);
  \draw (c) to node [above] {$e''$} (b);
\end{tikzpicture}
}
\end{center}
%\centerline{\mpfile{metric-res}{2}\hskip2cm
%\raisebox{2mm}{\mpfile{metric-res}{3}}}
  \caption{Parallel and series connection.
  All edges are directed from left to right.} % WE HAVE TO ADD ARROWS.}
  \label{pic2}
\end{figure}

% \label{s4}
Let us consider two simplest two-pole circuits
given in Figure \ref{pic2}.

%%%%\begin{figure}[h]
%%%%\centerline{\mpfile{metric-res}{2}\hskip2cm
%%%%\raisebox{2mm}{\mpfile{metric-res}{3}}}
%%%%  \caption{Parallel and series connection}\label{pic2}
%%%%\end{figure}

\begin{proposition}
\label{l-ps}
The resistances of these two circuits can be
determined, respectively, from formulas

\begin{equation}
\label{eq-ps}
\mu_{a,b}^{-s} = (\mu_{e'}^{-s} + \mu_{e''}^{-s})
\;\;\; \mbox{and} \;\;\;
\mu_{a,b}^{s/r} = (\mu_{e'}^{s/r} + \mu_{e''}^{s/r}).
\end{equation}
%\noindent
%respectively, where $e'$ and  $e''$  denote two edges of each circuit.
\end{proposition}

\proof
If  $r = s = 1$  then (\ref{eq-ps}) turns into
familiar high-school formulas.
The general case is just a little more difficult.
Without loss of generality let us assume that
$y_{a,b} = x_a - x_b \geq 0$.

\smallskip

In case of the parallel connection we obtain the following chain of equalities.

$$y^*_{a,b} = f_{a,b}(y_{a, b}) = \frac {y_{a,b}^r}{\mu_{a,b}^s} =
f_{e'}(y_{a,b}) + f_{e''}(y_{a,b}) =
\frac {y_{e'}^r} {\mu_{e'}^s} + \frac {y_{e''}^r} {\mu_{e''}^s} =
\frac {y_{a,b}^r} {\mu_{e'}^s} + \frac {y_{a,b}^r} {\mu_{e''}^s}.$$

Let us compare the third and the last terms;
dividing both by the numerator $y_{a,b}^r$  we arrive at (\ref{eq-ps}).

\medskip

In case of the series connection, let us start with determining
$x_c$  from the first Kirchhoff law:

$$y^*_{a,b} = f_{a,b}(y_{a, b}) =
\frac {y_{a,b}^r}{\mu_{a,b}^s} = \frac {(x_a - x_b)^r}{\mu_{a,b}^s}= $$

$$y^*_{e'} = f_{e'}(y_{e'}) =  f_{e'}(x_a - x_c) = \frac {(x_a - x_c)^r} {\mu_{e'}^s} =
y^*_{e''} = f_{e''}(y_{e''}) =  f_{e''}(x_c - x_b) = \frac {(x_c - x_b)^r} {\mu_{e''}^s}.$$

It is sufficient to compare the last and eighth terms to get

$$x_c = \frac{x_b \mu_{e'}^{s/r} + x_a \mu_{e''}^{s/r}}
{\mu_{e'}^{s/r} + \mu_{e''}^{s/r}}.$$

Then, let compare the last and forth terms, substitute
the obtained $x_c$, and get (\ref{eq-ps}).
\qed

\medskip

Now, let us consider the convolution  $\mu(t) = (\mu_{e'}^t +  \mu_{e''}^t)^{1/t}$;
it is well known and easy to see that

\begin{equation}
\label{convolution}
\mu(t) \rightarrow \max(\mu_{e'}, \mu_{e''}),
\;\mbox{as}\;
t \rightarrow + \infty,
\;\;\;\mbox{and}\;\;\;
\mu(t) \rightarrow \min(\mu_{e'}, \mu_{e''}),
\;\mbox{as}\;
t \rightarrow - \infty.
\end{equation}

\subsection*{Main four examples of resistance distances}
\label{s5}

Let us fix a weighted digraph   $G = (V, E, \mu)$  and
two strictly positive real parameters $r$  and  $s$.
As we proved, the obtained circuit
can be viewed as a quasi-metric space in which the distance
from   $a$  to  $b$  is defined as the effective resistance  $\mu_{a,b}$.
As announced in the introduction, this model
results in several interesting examples of quasi-metric
and quasi-ultrametric spaces.
Yet, to arrive to them we should allow for
$r$  and  $s/r$  to take values  $0$  and  $+\infty$.
More accurately, let  $r = r(t)$  and  $s = s(t)$  depend
on a real positive  parameter  $t$, 
or in other words, these two functions define a curve
in the positive quadrant  $s \geq 0, r \geq 0$.

% By Theorem \ref\ref{to},
We proved that resistances  $\mu_{a, b} (t)$  are well-defined
for every two nodes  $a, b \in V$  and each $t$.
We will show that, for the four limit transitions
listed below, limits
$\mu_{a, b} (t) = \lim_{t \rightarrow \infty} \mu_{a, b}(t)$,
exist for all  $a, b \in V$  and can be interpreted as follows:

\medskip

{\bf Example 1: the effective Ohm resistance of an electrical circuit.}

Let a weighted digraph  $G = (V, E, \mu)$  model an electrical circuit
in which  $\mu_e$  is the resistance of a directed edge
(semiconductor) $e$  and  $\; r(t) = s(t) \equiv 1$, or
more generally,  $\; r(t) \rightarrow 1$  and  $s(t) \rightarrow 1$, as  $t \rightarrow + \infty$.
Then,  $\mu_{a,b}$  is the effective Ohm resistance from  $a$  to  $b$.
For parallel and series connection of two directed edges
$e'$  and  $e''$, as in Figure~\ref{pic2},  we obtain, respectively,
$\;\; \mu_{a,b}^{-1} = \mu_{e'}^{-1} + \mu_{e''}^{-1}\; $  and
$\;\; \mu_{a,b} = \mu_{e'} + \mu_{e''},$ which is known from the high school.

\medskip

{\bf Example 2: the length of a shortest route.}

Let a weighted digraph  $G = (V, E, \mu)$  model a road network in which
$\mu_e$  is the length
(milage, traveling time, or gas consumption) of a one-way road  $e$.
Then,  $\mu_{a,b}$  can be viewed as the distance from
$a$  to  $b$, that is, the length of
a shortest directed path between them.
In this case, for parallel and series connection of
$e'$  and  $e''$,  we obtain, respectively,
$\mu_{a,b} = \min(\mu_{e'},\mu_{e''})\;\;$  and
$\;\; \mu_{a,b} = \mu_{e'} + \mu_{e''}.$
Hence, by (\ref{convolution}),  $\; -s(t) \rightarrow -\infty$  and  $s(t) \equiv r(t)$
for all  $t$, as in Figure \ref{pic0};
or more generally, $\; s(t) \rightarrow \infty\;$  and  $\; s(t)/r(t) \rightarrow 1$,
as  $t \rightarrow + \infty$.

\medskip

{\bf Example 3: the inverse width of a bottleneck route.}

Now, let digraph  $G = (V, E, \mu)$  model a system of one-way passages
(rivers, canals, bridges, etc.), where the conductance
$\lambda_e = \mu_e^{-1}$  is
the {\em "width"} of a passage  $e$, that is, the maximum size
(or tonnage) of a ship or a car that can pass  $e$, yet.
Then,  the effective conductance
$\lambda_{a,b} = \mu_{a,b}^{-1}$  is interpreted as the maximum width
of a (bottleneck) path between  $a$  and  $\;b$, that
is, the maximum size (or tonnage) of a ship or a car that can
still pass between terminals  $a$  and  $b$.
In this case,
$\;\; \lambda_{a,b} = \max(\lambda_{e'},\lambda_{e''})\;\;$
for the parallel connection and
$\;\; \lambda_{a,b} = \min(\lambda_{e'}, \lambda_{e''})$
for the series connection. Hence,
$s(t) \rightarrow \infty$  and   $\; s(t)/r(t) \rightarrow \infty$, as  $t \rightarrow \infty$;
in particular, $r$  might be bounded by a constant, $r(t) \leq const$, or just
$r(t) \equiv 1$  for all  $t$, as in Figure \ref{pic0}.

\medskip

{\bf Example 4: the inverse value of a maximal flow.}

Finally, let digraph  $G = (V, E, \mu)$  model
a pipeline or transportation network in which
the conductance $\lambda_e = \mu_e^{-1}$  is
the capacity of a one-way pipe or road  $e$.
Then,  $\lambda_{a,b} = \mu_{a,b}^{-1}$  is the capacity
of the whole two-pole network $(G, a, b)$  from terminal  $a$  to  $\;b$.
(Standardly, the capacity is defined as the amount of material
that can be transported through  $e$,
or from   $a$  to  $b$  in the whole circuit, per unit time.)
In this case,
$\;\; \lambda_{a,b} = \lambda_{e'} + \lambda_{e''}\;\;$
for the parallel connection and
$\;\; \lambda_{a,b} = \min(\lambda_{e'}, \lambda_{e''})$
for the series connection. Hence,
$-s(t) \equiv -1$  and  $s(t)/r(t) \rightarrow \infty$, that
is, $s(t) \equiv 1$  and  $r(t) \rightarrow 0$, as
in Figure \ref{pic0}, or
more generally, $\; s(t) \rightarrow 1$, while
$r(t) \rightarrow 0$, as  $t \rightarrow \infty$.

\begin{theorem}
\label{t2}
In all four examples, the limits
$\mu_{a, b} = \lim_{t \rightarrow +\infty} \mu_{a, b} (t) $  exist
and equal the corresponding distances from  $a$  to   $b$
for all $a, b \in V$.
In all four cases these distances define quasi-metric spaces
and in the last two - quasi-ultrametric spaces.
\end{theorem}

%--------------------------------------------------------------------------
%2. time
%$r(t) \rightarrow \infty$  and   $\; s(t)/r(t) \rightarrow 1$, as  $t \rightarrow \infty$;
%3. bottleneck
%$s(t) \rightarrow \infty$  and   $s(t)/r(t) \rightarrow \infty$, as  $t \rightarrow \infty$;
%4. flow
%$s(t) \rightarrow 1$, while $r(t) \rightarrow 0$, as  $t \rightarrow \infty$.
%---------------------------------------------------------------------------

\proof (sketch)
For Example 1 there is nothing to prove.
Also, for the series-parallel circuits
the statement is obvious in all cases.
It remans to consider Examples 2, 3,  and 4
for general circuits.
In each case our analysis will be based on
the minimum dissipation principle
(\ref{dissipation-e,dissipation}).

\medskip

For simplicity, we will omit argument  $t$  in
$r, s, \mu_e, \mu{a,b}, \lambda_e, \lambda_{a,b}, y^*_e, y_{a,b}$,
remembering, however, that all these variables
depend on  $t$, as indicated in the definitions  of Examples 2,3 and 4.

\bigskip

{\bf Example 2}.
In this case   $F_e(y^*_e) \sim  \mu_e  y^*_e$,
as  $t \rightarrow + \infty$ .

Hence,  % It is easy to see that
``moving some current to a shorter directed path"
reduces the total dissipation $F^*(y^*)$  when  $t$  is large enough.
More precisely, let
$p' = p'(s,t)$  and  $p''(s,t)$
be two directed paths  from  $s$  to  $t$  in  $G$  such that 
the first one is shorter, that is,  
\smallskip

$\mu(p') = \sum_{e \in p'} \mu_e < \sum_{e \in p''} \mu_e = \mu(p'')$.

\smallskip

Suppose that 
$\min (y^*_e \mid e \in p' \cup p'') = y^*_0 > 0$.
Obviously, dissipation  $F^*(y^*)$  will we reduced
by about  $(\mu(p'')  - \mu(p'))  y^*_0$
if we move  the current  $y^*_0$  from  $p''$  to  $p'$, that is,
we subtract  $y^*_0$  from  $y^*_e$  for each  $e \in p''$
and add it to  $y^*_e$  for each  $e \in p'$.

Recall that for every fixed  $t > 0$
each solution of the circuit  $(G,a,b)$
has a unique distribution of currents  $y^{*}$,
which minimizes the total dissipation  $F^*(y^*)$.
Thus, the above observation implies that
all currents will tend to the shortest
directed paths from  $a$  to  $b$,  as $t \rightarrow + \infty$.
Moreover,  $y_e^*$  well become just $0$  for
any  $e$  that does not belong to some shortest
$(a-b)$-path if  $t$  is large enough.

Clearly, by an arbitrary small perturbation of  $\mu_e$,
one can make the lengths of all
directed paths from  $a$  to $b$  distinct.
After such perturbation,
the shortest $(a-b)$-path  $p_0$  in $G$  becomes unique
and all currents outside of  it  become $0$,
that is,  $y^*_e = 0$  whenever  $e \not\in p_0$,
in particular, $\mu_{a,b}$  becomes  the length of  $p_0$
when  $t$  is large  enough.

\bigskip

{\bf Example 3}.
In this case  $F_e(y^*_e) \sim  \frac{1}{2} \mu_e^s  y^{*2}_e$,
where $s = s(t) \rightarrow + \infty$, as  $t \rightarrow + \infty$.
Thus, all currents will tend to the widest bottleneck
directed paths from  $a$  to  $b$,  as $t \rightarrow + \infty$.
Moreover,  $y_e^*$  well become just $0$  for
any  $e$  that does not belong to such a path if  $t$  is large enough.

Recall that
the {\em widest bottleneck}
directed  $(a-b)$-path  $p$  in  $G$
is defined as one maximizing
$\min(\mu_e \mid e \in p)$.

Unlike the shortest  $(a-b)$-path
the widest bottleneck one is
"typically"  not unique in  $G$.
Let us refine slightly this concept and
introduce the {\em lexicographically widest
bottleneck} directed $(a,b)$-path in  $G$.
To do so, consider all
widest bottleneck directed  $(a-b)$-paths in  $G$.
Among them choose those that maximize
the second smallest width, etc.
In several steps
(at most  $|V|$)  we will obtain the required path.

Clearly, by an arbitrary small perturbation of  $\mu_e$,
one can make the widths of all
directed paths from  $a$  to $b$  distinct.
Under this condition, the lexicographically widest
bottleneck directed $(a,b)$-path  $p_0$  in  $G$  becomes unique,
and all currents outside of it become $0$,
that is,  $y^*_e = 0$  whenever  $e \not\in p_0$,
in particular, $\mu_{a,b}$  becomes  the width  of  $p_0$
when  $t$  is large  enough.

\bigskip

{\bf Example 4}.
In this case, setting  $s =1$,  we obtain

\smallskip

$F_e(y^*_e) \sim  r y_e^* (\mu_e  y_e^*)^\frac{1}{r}  =
r y_e^* (\frac{y_e^*}{\lambda_e})^\frac{1}{r}$,
where  $r = r(t) \rightarrow 0$  as  $t \rightarrow +\infty$.

\smallskip

Let us recall the concept of the so-called {\em balanced flow} 
introduced in \cite{GG84} for the multi-pole circuits.
Given a weighted digraph  $G = (V, E)$   
consider the following boundary conditions: 
$x^*_v = x^{*0}_v$  for all $v \in V$. 

A flow  $y^*$  is called {\em satisfactory} if 
it satisfies these conditions. 
We will assume that such a flow exists. 
In particular, $\sum_{v \in V} x^{*0} = 0$  must hold. 

Two-pole boundary conditions (\ref{boundary*}),  
considered in the  present paper, form  
a special case of the multi-pole conditions.   
In this case a satisfactory flow exists  
if and only if digraph  $G$  
contains a directed path from  $a$ to $b$.  

Introduce resistances  $\mu_e$  
(and conductances  $\lambda_e = \mu_e^{-1}$)   
for all directed edges   $e  \in E$, 
thus getting a weighted multi-pole circuit $G = (V, E, \mu)$. 

Among all satisfactory flows choose all that
minimize  $\max (\frac{y_e^*}{\lambda_e} \mid  e \in E)$;
then among them choose all that minimize
the second largest value of  $\frac{y_e^*}{\lambda_e}$, etc.
For all $e \in E$  order the ratios $\frac{y_e^*}{\lambda_e}$  non-decreasingly.
A satisfactory flow  $y^*$  realizing the lex-min over all such vectors
is called {\em balanced}.

Let us briefly recall the 
% linear time 
algorithm from \cite{GG84} 
constructing a balanced flow in a multi-pole circuit.
A cut in a digraph  $G = (V, E)$,
is defined as an ordered pair  $C = (V', V")$  
that partitions  $V$  properly, that is, 

\smallskip 

$V' \neq \emptyset, V'' \neq \emptyset, V' \cap V'' = \emptyset$  and  $V' \cup V'' = V$. 

\smallskip 

We say that a directed edge  $e = (u,w)$  is 
in  $C$  if $u \in V'$  and  $w \in V''$.  
Consider a multi-pole circuit  
defined by a weighted digraph  $G = (V, E, \mu)$  and 
boundary conditions  $x^*_v = x^{*0}_v, \; v \in V$; 
fix a cut  $C$  in  $G$. 
Its  {\em  deficiency} and  {\em capacity} are  defined by formulas:   

\smallskip 

$D(C) = \sum_{v \in V'} x^{*0} = - \sum_{v \in V''} x^{*0}; \;\; 
\Lambda(C) = \sum_{e \in C} \lambda_e$. 

\smallskip 

Note that  $\Lambda(C) \geq 0$,  by this definition, while 
$D(C)$  may be negative. 

A multi-pole problem has no  satisfactory vector if and only if 
there exists a cut  $C$  such that  $R(C) = +\infty$, or in other words, 
such that  $D(C) > 0$  and  $\Lambda(C) = 0$.
In a two-pole circuit 
this happens  if and only if 
there is  no directed path from  $a$  to  $b$. 
Note also that in the two-pole case we have
$D(C) = 0$  whenever  $a,b \in V'$  or  $a,b \in V''$.

Cut  $C$  is called {\em critical} if
it realizes the maximum of the ratio 
$R(C) = \frac{D(C)}{\lambda(C)}$. 
%  among all cuts of  $G$.
% A critical cut can be found in linear time.  
Choose such  $C$  and set  
$y_e^* = R(C) \lambda_e$   for each  $e  \in C$. 
Then, $\frac{y_e^*}{\lambda_e}$
take  the same value  $R(C)$   for all $e \in C$
and we have  $\sum_{e \in C} y_e^* = \lambda(C)$.

Reduce digraph  $G$  eliminating all edges of  $C$  from it.  
Recompute new multi-pole boundary conditions 
in the obtained reduced digraph  $G'$  
taking into account flows  $y_e^*$  on the deleted edges  $e \in C$.   
Then, find a critical cut  $C'$  in   $G'$, and repeat. 
It is shown in \cite{GG84} that the values 
$R, R',  \ldots$  are monotone non-increasing.     
Hence, the ratio  $\frac{y_e^*}{\lambda_e}$   
takes the  largest values on  $e \in C$  at the first stage.  
It is also shown in \cite{GG84} that the  balanced flow is unique.  
%  and  found by the above algorithm in polynomial time.  

\medskip 

Recall that  
$F_e(y^*_e) \sim  (r y_e^*) (\frac{y_e^*}{\lambda_e})^\frac{1}{r}$,
where $r = r(t) \rightarrow 0$  as  $t \rightarrow +\infty$.
From this we conclude that  
a satisfactory flow  $y^*$, 
minimizing the total dissipation $F^*(y^*)$,  
becomes a balanced flow when  $t$  is large enough.

\smallskip

A satisfactory flow  $y^*$  is called {\em feasible} 
if  $y_e^* \leq \lambda_e$.
% (Recall that  $y^*_e \geq 0$ for all $e \in E$ 
% and we assume that a satisfactory flow exists.) 
Obviously, the (unique) balanced flow  $y^{*0}$ is feasible 
whenever a latter exists. 
We saw that in this case   $y^{*0}$  is a solution 
of the multi-pole network problem, since  $y^{*0}$  
minimizes  $F^*(y^*)$   when $t$  is large enough.   
 
In particular, this is true for the two-pole problems.
In this  case let us set  $x_a^{*0} = - x_b^{*0} = \Lambda(a,b)$, 
where  $\Lambda(a,b)$  is the capacity 
from  $a$ to $b$  of the circuit  $G = (V, E, \mu)$   
with poles  $a$  and  $b$.   
Obviously, on the first stage of the algorithm we obtain 
a cut  $C$  of the unit ratio $R(C) = 1$, that is,  
$y^{*0}_e = \lambda e$  for all  $e \in C$.  
Thus, $\lambda(a,b) \rightarrow \Lambda(a,b)$,  as $t \rightarrow +\infty$ 
\qed

\bigskip

\section*{Acknowledgements}
The paper was prepared within the framework
of the HSE University Basic Research Program and
funded by the RSF grant  20-11-20203;
The author is thankful to Endre Boros and
Mert G\"{u}rb\"{u}zbalaban for many helpful remarks.

\end{document}